\documentclass{elsarticle}
\usepackage[T1]{fontenc}
\usepackage[latin9]{inputenc}
\usepackage{prettyref}
\usepackage{enumitem}
\usepackage{amsthm}
\usepackage{amsmath}
\usepackage{amssymb}
\usepackage[all]{xy}

\makeatletter
\theoremstyle{plain}
\newtheorem{thm}{\protect\theoremname}
  \theoremstyle{plain}
  \newtheorem{prop}[thm]{\protect\propositionname}
 \newlist{casenv}{enumerate}{4}
 \setlist[casenv]{leftmargin=*,align=left,widest={iiii}}
 \setlist[casenv,1]{label={{\itshape\ \casename} \arabic*.},ref=\arabic*}
 \setlist[casenv,2]{label={{\itshape\ \casename} \roman*.},ref=\roman*}
 \setlist[casenv,3]{label={{\itshape\ \casename\ \alph*.}},ref=\alph*}
 \setlist[casenv,4]{label={{\itshape\ \casename} \arabic*.},ref=\arabic*}
  \theoremstyle{plain}
  \newtheorem{lem}[thm]{\protect\lemmaname}
  \theoremstyle{plain}
  \newtheorem{cor}[thm]{\protect\corollaryname}
  \theoremstyle{definition}
  \newtheorem{example}[thm]{\protect\examplename}

\usepackage{braket}

\xyoption{cmtip}
\SelectTips{cm}{}

\newrefformat{prop}{Proposition \ref{#1}}

\makeatother

  \providecommand{\corollaryname}{Corollary}
  \providecommand{\examplename}{Example}
  \providecommand{\lemmaname}{Lemma}
  \providecommand{\propositionname}{Proposition}
 \providecommand{\casename}{Case}
\providecommand{\theoremname}{Theorem}

\begin{document}

\begin{frontmatter}{}

\title{Nonstandard homology theory for uniform spaces}

\author{Takuma Imamura}

\address{Department of Mathematics\\
University of Toyama\\
3190 Gofuku, Toyama 930-8555, Japan}

\ead{s1240008@ems.u-toyama.ac.jp}
\begin{abstract}
We introduce a new homology theory of uniform spaces, provisionally
called $\mu$-homology theory. Our homology theory is based on hyperfinite
chains of microsimplices. This idea is due to McCord. We prove that
$\mu$-homology theory satisfies the Eilenberg-Steenrod axioms. The
characterization of chain-connectedness in terms of $\mu$-homology
is provided. We also introduce the notion of S-homotopy, which is
weaker than uniform homotopy. We prove that $\mu$-homology theory
satisfies the S-homotopy axiom, and that every uniform space can be
S-deformation retracted to a dense subset. It follows that for every
uniform space $X$ and any dense subset $A$ of $X$, $X$ and $A$
have the same $\mu$-homology. We briefly discuss the difference and
similarity between $\mu$-homology and McCord homology.\end{abstract}

\end{frontmatter}{}

\section{Introduction}

McCord \cite{McC72} developed a homology of topological spaces using
nonstandard methods. McCord's theory is based on hyperfinite chains
of microsimplices. Intuitively, microsimplices are abstract simplices
with infinitesimal diameters. Garavaglia \cite{Gar78} proved that
McCord homology coincides with \v{C}ech homology for compact spaces.
\v{Z}ivaljevi\'c \cite{Ziv87} proved that McCord cohomology also coincides
with \v{C}ech cohomology for locally contractible paracompact spaces.
Korppi \cite{Kor10} proved that McCord homology coincides with \v{C}ech
homology with compact supports for regular Hausdorff spaces.

In this paper, we introduce a new microsimplicial homology theory
of uniform spaces, provisionally called $\mu$-homology theory. $\mu$-homology
theory satisfies the Eilenberg-Steenrod axioms. Vanishing of the $0$-th
reduced $\mu$-homology characterizes chain-connectedness. We also
introduce the notion of S-homotopy, which is weaker than uniform homotopy.
$\mu$-homology theory satisfies the S-homotopy axiom. Hence $\mu$-homology
is an S-homotopy invariant. Every uniform space can be S-deformation
retracted to a dense subset. It follows that for every uniform space
$X$ and any dense subset $A$ of $X$, $X$ and $A$ have the same
$\mu$-homology. We briefly discuss the difference and similarity
between $\mu$-homology and McCord homology.

\section{Preliminaries}

The basics of nonstandard analysis are assumed. We fix a universe
$\mathbb{U}$, the standard universe, satisfying sufficiently many
axioms of ZFC. All standard objects we consider belong to $\mathbb{U}$.
We also fix an elementary extension $^{\ast}\mathbb{U}$ of $\mathbb{U}$,
the internal universe, that is $\left|\mathbb{U}\right|^{+}$-saturated.
The map $x\mapsto{}^{\ast}x$ denotes the elementary embedding from
$\mathbb{U}$ into $^{\ast}\mathbb{U}$. We say ``by transfer'' to
indicate the use of the elementary equivalence between $\mathbb{U}$
and $^{\ast}\mathbb{U}$. We say ``by saturation'' when using the
saturation property of $^{\ast}\mathbb{U}$.

Let us enumerate some well-known facts of nonstandard topology. Let
$X$ be a topological space. The monad of $x\in X$ is $\mu\left(x\right)=\bigcap\set{^{\ast}U|x\in U\in\tau}$,
where $\tau$ is the topology of $X$. A subset $U$ of $X$ is open
if and only if $\mu\left(x\right)\subseteq{}^{\ast}U$ for all $x\in U$.
A subset $F$ of $X$ is closed if and only if $\mu\left(x\right)\cap{}^{\ast}F\neq\varnothing$
implies $x\in F$ for all $x\in X$. A subset $K$ of $X$ is compact
if and only if for any $x\in{}^{\ast}K$ there is a $y\in K$ with
$x\in\mu\left(y\right)$. A map $f:X\to Y$ of topological spaces
is continuous at $x\in X$ if and only if for any $y\in\mu\left(x\right)$
we have $^{\ast}f\left(y\right)\in\mu\left(f\left(x\right)\right)$.

Next, let $X$ be a uniform space. Two points $x,y$ of $^{\ast}X$
are said to be infinitely close, denoted by $x\approx y$, if for
any entourage $U$ of $X$ we have $\left(x,y\right)\in{}^{\ast}U$.
$\approx$ is an equivalence relation on $^{\ast}X$. The monad of
$x\in X$ is equal to $\mu\left(x\right)=\set{y\in{}^{\ast}X|x\approx y}$.
Thus, in the case of uniform spaces, one can define the monad of $x\in{}^{\ast}X$.
For each entourage $U$ of $X$, the $U$-neighbourhood of $x\in X$
is $U\left[x\right]=\set{y\in X|\left(x,y\right)\in U}$. A map $f:X\to Y$
of uniform spaces is uniformly continuous if and only if $x\approx y$
implies $^{\ast}f\left(x\right)\approx{}^{\ast}f\left(y\right)$ for
all $x,y\in{}^{\ast}X$.

Let $\set{X_{i}}_{i\in I}$ be a family of uniform spaces. Let $P$
be the product $\prod_{i\in I}X_{i}$ of $\set{X_{i}}_{i\in I}$,
and let $Q$ be the coproduct $\coprod_{i\in I}X_{i}$ of $\set{X_{i}}_{i\in I}$.
Let $\approx_{X}$ denote the ``infinitely close'' relation of a uniform
space $X$. For any $x,y\in P$, $x\approx_{P}y$ if and only if $x\left(i\right)\approx_{X_{i}}y\left(i\right)$
for all $i\in I$. For any $x,y\in Q$, $x\approx_{Q}y$ if and only
if there is an $i\in I$ such that $x,y\in X_{i}$ and $x\approx_{X_{i}}y$.

\section{Definition of $\mu$-homology theory}

Let $X$ be a uniform space and $G$ an internal abelian group. We
denote by $C_{p}X$ the internal free abelian group generated by $^{\ast}X^{p+1}$,
and by $C_{p}\left(X;G\right)$ the internal abelian group of all
internal homomorphisms from $C_{p}X$ to $G$. Each member of $C_{p}\left(X;G\right)$
can be represented in the form $\sum_{i=0}^{n}g_{i}\sigma_{i}$, where
$\set{g_{i}}_{i=0}^{n}$ is an internal hyperfinite sequence of members
of $G$, and $\set{\sigma_{i}}_{i=0}^{n}$ is an internal hyperfinite
sequence of members of $^{\ast}X^{p+1}$. A member $\left(a_{0},\ldots,a_{p}\right)$
of $^{\ast}X^{p+1}$ is called a \emph{$p$-microsimplex} if $a_{i}\approx a_{j}$
for all $i,j\leq p$, or equivalently, $\mu\left(a_{0}\right)\cap\cdots\cap\mu\left(a_{p}\right)\neq\varnothing$.
A member of $C_{p}\left(X;G\right)$ is called a \emph{$p$-microchain}
if it can be represented in the form $\sum_{i=0}^{n}g_{i}\sigma_{i}$,
where $\set{g_{i}}_{i=0}^{n}$ is an internal hyperfinite sequence
of members of $G$, and $\set{\sigma_{i}}_{i=0}^{n}$ is an internal
hyperfinite sequence of $p$-microsimplices. We denote by $M_{p}\left(X;G\right)$
the subgroup of $C_{p}\left(X;G\right)$ consisting of all $p$-microchains.
The boundary map $\partial_{p}:M_{p}\left(X;G\right)\to M_{p-1}\left(X;G\right)$
is defined by
\[
\partial_{p}\left(a_{0},\ldots,a_{p}\right)=\sum_{i=0}^{p}\left(-1\right)^{i}\left(a_{0},\ldots,\hat{a}_{i},\ldots,a_{p}\right).
\]
More precisely, we first define an internal map $\partial_{p}':C_{p}\left(X;G\right)\to C_{p-1}\left(X;G\right)$
by the same equation. We see that $\partial_{p}'\left(M_{p}\left(X;G\right)\right)\subseteq M_{p-1}\left(X;G\right)$.
$\partial_{p}$ is defined by the restriction of $\partial_{p}'$
to $M_{p}\left(X;G\right)$. Thus $M_{\bullet}\left(X;G\right)$ forms
a chain complex. 

Let $f:X\to Y$ be a uniformly continuous map. By the nonstandard
characterization of uniform continuity, we see that for every $p$-microsimplex
$\left(a_{0},\ldots,a_{p}\right)$ on $X$, $\left(^{\ast}f\left(a_{0}\right),\ldots{}^{\ast}f\left(a_{p}\right)\right)$
is a $p$-microsimplex on $Y$. The induced homomorphism $M_{\bullet}\left(f;G\right):M_{\bullet}\left(X;G\right)\to M_{\bullet}\left(Y;G\right)$
of $f$ is defined by
\[
M_{p}\left(f;G\right)\left(a_{0},\ldots,a_{p}\right)=\left(^{\ast}f\left(a_{0}\right),\ldots{}^{\ast}f\left(a_{p}\right)\right).
\]
Thus we have the functor $M_{\bullet}\left(\cdot;G\right)$ from the
category of uniform spaces to the category of chain complexes. \emph{$\mu$-homology
theory} is the composition of functors $H_{\bullet}\left(\cdot;G\right)=H_{\bullet}M_{\bullet}\left(\cdot;G\right)$,
where $H_{\bullet}$ in the right hand side is the ordinary homology
theory of chain complexes.

Let $X$ be a uniform space and $A$ be a subset of $X$. The induced
homomorphism $M_{p}\left(i;G\right):M_{p}\left(A;G\right)\to M_{p}\left(X;G\right)$
of the inclusion map $i:A\hookrightarrow X$ is injective. Let us
identify $M_{\bullet}\left(A;G\right)$ with a subchain complex of
$M_{\bullet}\left(X;G\right)$ and define
\[
M_{\bullet}\left(X,A;G\right)=\frac{M_{\bullet}\left(X;G\right)}{M_{\bullet}\left(A;G\right)}.
\]
Every uniformly continuous map $f:\left(X,A\right)\to\left(Y,B\right)$
induces a homomorphism $M_{\bullet}\left(f;G\right):\left(X,A;G\right)\to\left(Y,B;G\right)$.
Thus $M_{\bullet}\left(\cdot,\cdot;G\right)$ is a functor from the
category of pairs of uniform spaces to the category of chain complexes.
\emph{Relative $\mu$-homology theory} is the composition of functors
$H_{\bullet}\left(\cdot,\cdot;G\right)=H_{\bullet}M_{\bullet}\left(\cdot,\cdot;G\right)$.

\section{Eilenberg-Steenrod axioms}

In this section, we will verify that $\mu$-homology theory satisfies
the Eilenberg-Steenrod axioms: uniform homotopy, exactness, weak excision,
dimension, and finite additivity.

Recall that two uniformly continuous maps $f,g:X\to Y$ are said to
be \emph{uniformly homotopic} if there is a uniformly continuous map
$h:X\times\left[0,1\right]\to Y$, called a \emph{uniform homotopy}
between $f$ and $g$, such that $h\left(\cdot,0\right)=f$ and $h\left(\cdot,1\right)=g$.
\begin{prop}[Uniform homotopy]
\label{prop:uniform-homotopy-axiom}If two uniformly continuous maps
$f,g:X\to Y$ are uniformly homotopic, then the induced homomorphisms
$M_{\bullet}\left(f;G\right)$ and $M_{\bullet}\left(g;G\right)$
are chain homotopic. Hence $H_{\bullet}\left(f;G\right)=H_{\bullet}\left(g;G\right)$.
This also holds for relative $\mu$-homology.\end{prop}
\begin{proof}
Let $h$ be a uniform homotopy between $f$ and $g$. Fix an infinite
hypernatural number $N$. Define $h_{i}={}^{\ast}h\left(\cdot,i/N\right)$.
For each hypernatural number $i\leq N$, we define a map $P_{i,p}:M_{p}\left(X;G\right)\to M_{p+1}\left(Y;G\right)$
by letting
\[
P_{i,p}\left(a_{0},\ldots,a_{p}\right)=\sum_{j=0}^{p}\left(-1\right)^{j}\left(h_{i}\left(a_{0}\right),\ldots,h_{i}\left(a_{j}\right),h_{i+1}\left(a_{j}\right),\ldots,h_{i+1}\left(a_{p}\right)\right).
\]
Note that the hyperfinite sequence $\set{P_{i,p}u}_{i=0}^{N-1}$ is
internal. Hence the hyperfinite sum $P_{p}u=\sum_{i=0}^{N-1}P_{i,p}u$
exists. Thus we obtain the prism map $P_{p}:M_{p}\left(X;G\right)\to M_{p+1}\left(Y;G\right)$.
Let us verify that $P_{\bullet}$ is a chain homotopy between $M_{\bullet}\left(f;G\right)$
and $M_{\bullet}\left(g;G\right)$.
\begin{align*}
 & \partial_{p+1}P_{i,p}\left(a_{0},\ldots,a_{p}\right)=\\
 & \qquad\sum_{k<j}\left(-1\right)^{j+k}\left(h_{i}\left(a_{0}\right),\ldots,\widehat{h_{i}\left(a_{k}\right)},\ldots,h_{i}\left(a_{j}\right),h_{i+1}\left(a_{j}\right),\ldots,h_{i+1}\left(a_{p}\right)\right)\\
 & \qquad-\sum_{j<k}\left(-1\right)^{j+k}\left(h_{i}\left(a_{0}\right),\ldots,h_{i}\left(a_{j}\right),h_{i+1}\left(a_{j}\right),\ldots,\widehat{h_{i+1}\left(a_{k}\right)},\ldots,h_{i+1}\left(a_{p}\right)\right)\\
 & \qquad+\left(h_{i+1}\left(a_{0}\right),\ldots,h_{i+1}\left(a_{p}\right)\right)-\left(h_{i}\left(a_{0}\right),\ldots,h_{i}\left(a_{p}\right)\right),\\
 & P_{i,p-1}\partial_{p}\left(a_{0},\ldots,a_{p}\right)=\\
 & \qquad\sum_{j<k}\left(-1\right)^{j+k}\left(h_{i}\left(a_{0}\right),\ldots,h_{i}\left(a_{j}\right),h_{i+1}\left(a_{j}\right),\ldots,\widehat{h_{i+1}\left(a_{k}\right)},\ldots,h_{i+1}\left(a_{p}\right)\right)\\
 & \qquad-\sum_{k<j}\left(-1\right)^{j+k}\left(h_{i}\left(a_{0}\right),\ldots,\widehat{h_{i}\left(a_{k}\right)},\ldots,h_{i}\left(a_{j}\right),h_{i+1}\left(a_{j}\right),\ldots,h_{i+1}\left(a_{p}\right)\right).
\end{align*}
Thus we obtain 
\[
\left(\partial_{p+1}P_{i,p}+P_{i,p-1}\partial_{p}\right)\left(a_{0},\ldots,a_{p}\right)=\left(h_{i+1}\left(a_{0}\right),\ldots,h_{i+1}\left(a_{p}\right)\right)-\left(h_{i}\left(a_{0}\right),\ldots,h_{i}\left(a_{p}\right)\right)
\]
and $\partial_{p+1}P_{p}+P_{p}\partial_{p}=M_{p}\left(g;G\right)-M_{p}\left(f;G\right)$.\end{proof}
\begin{prop}[Exactness]
Let $X$ be a uniform space and $A$ a subset of $X$. The sequence
\[
\xymatrix{0\ar[r] & M_{\bullet}\left(A;G\right)\ar[r]^{i_{\bullet}} & M_{\bullet}\left(X;G\right)\ar[r]^{j_{\bullet}} & M_{\bullet}\left(X,A;G\right)\ar[r] & 0}
\]
is exact, where $i_{\bullet}$ is the inclusion map and $j_{\bullet}$
is the projection map. Moreover, the above short exact sequence splits.\end{prop}
\begin{proof}
The first part is immediate from the definition. We will construct
a right inverse $s_{\bullet}$ of $j_{\bullet}$. Let $\upsilon\in M_{p}\left(X,A;G\right)$.
Choose a representative $u=\sum_{i}g_{i}\sigma_{i}\in M_{p}\left(X;G\right)$.
Since $\sum_{\sigma_{i}\subseteq{}^{\ast}A}g_{i}\sigma_{i}\in M_{p}\left(A;G\right)$,
$u'=u-\sum_{\sigma_{i}\subseteq{}^{\ast}A}g_{i}\sigma_{i}$ is also
a representative of $\upsilon$. $u'$ is uniquely determined by $\upsilon$
and does not depend on the choice of $u$. Define $s_{p}\left(\upsilon\right)=u'$.
It is easy to see that $s_{\bullet}$ is a right inverse of $j_{\bullet}$.\end{proof}
\begin{prop}[Weak excision]
Let $X$ be a uniform space. Let $A$ and $B$ be subsets of $X$
such that $X=\mathrm{int}A\cup\mathrm{int}B$. If either $A$ or $B$
is compact, then the inclusion map $i:\left(A,A\cap B\right)\hookrightarrow\left(X,B\right)$
induces the isomorphism $H_{\bullet}\left(i;G\right):H_{\bullet}\left(A,A\cap B;G\right)\cong H_{\bullet}\left(X,B;G\right)$.\end{prop}
\begin{proof}
It suffices to show the following two inclusions:
\begin{enumerate}
\item $M_{p}\left(A;G\right)\cap M_{p}\left(B;G\right)\subseteq M_{p}\left(A\cap B;G\right)$,
\item $M_{p}\left(X;G\right)\subseteq M_{p}\left(A;G\right)+M_{p}\left(B;G\right)$.
\end{enumerate}

\noindent The first inclusion is clear. We will only prove the second
inclusion. Suppose $u=\sum_{i}g_{i}\sigma_{i}\in M_{p}\left(X;G\right)$.
If each $\sigma_{i}$ is contained in either $^{\ast}A$ or $^{\ast}B$,
then $u\in M_{p}\left(A;G\right)+M_{p}\left(B;G\right)$.
\begin{casenv}
\item \noindent $A$ is compact. Suppose that $\sigma_{i}$ is not contained
in $^{\ast}B$. Then $\sigma_{i}$ intersects $^{\ast}A$. By the
nonstandard characterization of compactness, there is an $x\in A$
such that all vertices of $\sigma_{i}$ are infinitely close to $x$.
$x$ must belong to $\mathrm{int}A$. Otherwise, by the nonstandard
characterization of open sets, $\sigma_{i}\subseteq\mu\left(x\right)\subseteq{}^{\ast}B$,
a contradiction. Hence $\sigma_{i}\subseteq\mu\left(x\right)\subseteq{}^{\ast}A$.
\item $B$ is compact. Suppose that $\sigma_{i}$ is not contained in $^{\ast}A$.
Then $\sigma_{i}$ intersects $^{\ast}B$. There is an $x\in B$ such
that all vertices of $\sigma_{i}$ are infinitely close to $x$. $x$
must contained in $\mathrm{int}B$. Otherwise, $\sigma_{i}\subseteq\mu\left(x\right)\subseteq{}^{\ast}A$,
a contradiction. Hence $\sigma_{i}\subseteq\mu\left(x\right)\subseteq{}^{\ast}B$.
\end{casenv}
\end{proof}
\begin{prop}[Dimension]
If $X$ is the one-point space, then
\[
H_{p}\left(X;G\right)=\begin{cases}
G, & p=0,\\
0, & p\neq0.
\end{cases}
\]
\end{prop}
\begin{proof}
Immediate by definition.\end{proof}
\begin{prop}[Finite additivity]
\label{prop:finite-additivity-axiom}Let $X=\coprod_{i=0}^{n}X_{i}$
be a finite coproduct of uniform spaces. Then $M_{\bullet}\left(X;G\right)=\bigoplus_{i=0}^{n}M_{\bullet}\left(X_{i};G\right)$.
This also holds for relative $\mu$-homology.\end{prop}
\begin{proof}
Suppose $u=\sum_{i}g_{i}\sigma_{i}\in M_{p}\left(X;G\right)$. Any
two points in different components of $^{\ast}X$ are not infinitely
close. Each $\sigma_{i}$ is contained in one and only one of the
$^{\ast}A_{i}$. Hence $u\in\bigoplus_{i=0}^{n}M_{\bullet}\left(X_{i};G\right)$.
\end{proof}

\section{Chain-connectedness}

Consider the augmented chain complex $\tilde{M}_{\bullet}\left(X;G\right)$:
\[
\xymatrix{G & M_{0}\left(X;G\right)\ar[l]_{\varepsilon} & M_{1}\left(X;G\right)\ar[l]_{\partial_{1}} & M_{2}\left(X;G\right)\ar[l]_{\partial_{2}} & \cdots\cdots\ar[l]}
\]
where $\varepsilon$ is the augmentation map $\varepsilon\sum_{i}g_{i}\sigma_{i}=\sum_{i}g_{i}$.
\emph{Reduced $\mu$-homology theory} is defined by $\tilde{H}_{\bullet}\left(\cdot;G\right)=H_{\bullet}\tilde{M}_{\bullet}\left(\cdot;G\right)$.
For each $p>0$, $\tilde{H}_{p}\left(X;G\right)$ is identical to
$H_{p}\left(X;G\right)$. $\tilde{H}_{0}\left(X;G\right)$ is a subgroup
of $H_{0}\left(X;G\right)$.
\begin{prop}
\label{prop:ses-of-aug-cc}Let $X$ be a nonempty uniform space. The
sequence
\[
\xymatrix{0\ar[r] & \tilde{H}_{0}\left(X;G\right)\ar[r]^{i} & H_{0}\left(X;G\right)\ar[r]^{\bar{\varepsilon}} & G\ar[r] & 0}
\]
is exact and splits, where $i$ is the inclusion map and $\bar{\varepsilon}$
is the map induced by the augmentation map $\varepsilon$.\end{prop}
\begin{proof}
The well-definedness of $\bar{\varepsilon}$ follows from $\mathrm{im}\partial_{1}\subseteq\ker\varepsilon$.
Since $i$ is injective, the sequence is exact at $\tilde{H}_{0}\left(X;G\right)$.
The exactness at $H_{0}\left(X;G\right)$ is clear. Let us fix a point
$x$ of $X$. For any $g\in G$, we have $\bar{\varepsilon}\left[gx\right]=\varepsilon\left(gx\right)=g$.
Hence $\bar{\varepsilon}$ is surjective. The exactness at $G$ is
proved. Define a map $s:G\to H_{0}\left(X;G\right)$ by letting $s\left(g\right)=\left[gx\right]$.
Clearly $s$ is a right inverse of $\bar{\varepsilon}$.
\end{proof}
The $0$-th $\mu$-homology relates to the notion of chain-connectedness
of uniform spaces. Recall that a uniform space $X$ is said to be
\emph{chain-connected} if it is $U$-connected for all entourage $U$
of $X$. Here $X$ is said to be \emph{$U$-connected} if for any
$x,y\in X$ there is a finite sequence $\set{x_{i}}_{i=0}^{n}$ of
points of $X$, called a \emph{$U$-chain}, such that $x_{0}=x$,
$x_{n}=y$ and $\left(x_{i},x_{i+1}\right)\in U$ for all $i<n$.
A chain-connected space is also called a \emph{well-chained} space.
Every connected uniform space is chain-connected. The converse is
not true, e.g., the real line without one point $\mathbb{R}\setminus\set{0}$
is chain-connected but not connected. We can easily get the following
characterization.
\begin{lem}
\label{lem:chain-connectedness-lemma}A uniform space $X$ is chain-connected
if and only if for any $x,y\in{}^{\ast}X$ there is an internal hyperfinite
sequence $\set{x_{i}}_{i=0}^{n}$ of points of $^{\ast}X$, called
an infinitesimal chain, such that $x_{0}=x$, $x_{n}=y$ and $x_{i}\approx x_{i+1}$
for all $i<n$.\end{lem}
\begin{proof}
Suppose first that $X$ is chain-connected. By saturation, there is
an internal entourage $U$ of $^{\ast}X$ with $U\subseteq\left(\approx\right)$.
By transfer, for any $x,y\in{}^{\ast}X$ there is an internal $U$-chain
connecting $x$ to $y$. This is an infinitesimal chain.

Conversely, suppose that for any $x,y\in{}^{\ast}X$ there is an infinitesimal
chain connecting $x$ to $y$. Let $x,y$ be any points of $X$. Let
$U$ be any entourage of $X$. There is an infinitesimal chain connecting
$x$ to $y$. This is an internal $U$-chain connecting $x$ to $y$.
By transfer, there is a $U$-chain connecting $x$ to $y$.\end{proof}
\begin{thm}
\label{thm:0hom-of-cc-space-vanishes}If $X$ is chain-connected then
$\tilde{H}_{0}\left(X;G\right)=0$.\end{thm}
\begin{proof}
It suffices to show that the sequence
\[
\xymatrix{G & M_{0}\left(X;G\right)\ar[l]_{\varepsilon} & M_{1}\left(X;G\right)\ar[l]_{\partial_{1}}}
\]
is exact. The only nontrivial part is $\ker\varepsilon\subseteq\mathrm{im}\partial_{1}$.
Let $\sigma,\tau$ be $0$-microsimplices. We shall identify a point
$x$ of $^{\ast}X$ with the $0$-microsimplex $\left(x\right)$.
By \prettyref{lem:chain-connectedness-lemma}, there is an infinitesimal
chain $\set{x_{i}}_{i=0}^{n}$ connecting $\sigma$ to $\tau$. Since
$\set{\left(x_{i},x_{i+1}\right)}_{i=0}^{n-1}$ is an internal hyperfinite
sequence of $1$-microsimplices, the $1$-microchain $\sum_{i=0}^{n-1}\left(x_{i},x_{i+1}\right)$
is well-defined, and $\tau-\sigma=\partial_{1}\sum_{i=0}^{n-1}\left(x_{i},x_{i+1}\right)\in\mathrm{im}\partial_{1}$.

Suppose $u=\sum_{i=0}^{n}g_{i}\sigma_{i}\in\ker\varepsilon$. By induction,
we have
\[
u=\sum_{i=0}^{n}\sum_{j=0}^{i}g_{j}\left(\sigma_{i}-\sigma_{i+1}\right)+\sum_{j=0}^{n}g_{j}\sigma_{n}=\sum_{i=0}^{n}\sum_{j=0}^{i}g_{j}\left(\sigma_{i}-\sigma_{i+1}\right).
\]
Hence $u\in\mathrm{im}\partial_{1}$.
\end{proof}
The converse is also true.
\begin{thm}
\label{thm:if-0hom-vanishes-then-cc}Suppose that $G$ is nontrivial.
If $\tilde{H}_{0}\left(X;G\right)=0$ then $X$ is chain-connected.\end{thm}
\begin{proof}
The support of a $p$-microchain $u$ is $\set{\sigma\in{}^{\ast}X^{p+1}|u\left(\sigma\right)\neq0}$.
We denote by $V_{u}$ the set of all vertices of members of the support
of $u$. $V_{u}$ is an internal hyperfinite set. Let $u$ be a $1$-microchain.
We say that a vertex $x$ is accessible to a vertex $y$ on $u$,
denoted by $x\underset{u}{\leftrightarrow}y$, if there is an internal
hyperfinite sequence $\set{x_{i}}_{i=0}^{n}$ in $V_{u}$ such that
$x,y\in\set{x_{i}}_{i=0}^{n}$ and for each $i<n$ either $\left(x_{i},x_{i+1}\right)$
or $\left(x_{i+1},x_{i}\right)$ belongs to the support of $u$. Clearly
$\underset{u}{\leftrightarrow}$ is an internal equivalence relation
on $V_{u}$ and $\underset{u}{\leftrightarrow}$ implies $\approx$.

Let $x,y\in{}^{\ast}X$. We will show that $x$ and $y$ can be connected
by an infinitesimal chain. Since $\tilde{H}_{0}\left(X;G\right)=0$,
there is a $u=\sum_{i}g_{i}\sigma_{i}\in M_{1}\left(X;G\right)$ with
$\partial_{1}u=x-y$. Let $V_{1},\ldots,V_{n}$ be all equivalence
classes of $\underset{u}{\leftrightarrow}$. Suppose for a contradiction
that $x$ and $y$ belong to different equivalence classes $V_{\xi}$
and $V_{\eta}$, respectively. Then,
\begin{align*}
\partial_{1}u & =\sum_{i}g_{i}\partial_{1}\sigma_{i}\\
 & =\sum_{j=1}^{n}\sum_{\sigma_{i}\subseteq V_{j}}g_{i}\partial_{1}\sigma_{i}\\
 & =\sum_{\sigma_{i}\subseteq V_{\xi}}g_{i}\partial_{1}\sigma_{i}+\sum_{\sigma_{i}\subseteq V_{\eta}}g_{i}\partial_{1}\sigma_{i}+\sum_{\substack{j=1\\
j\neq\xi,\eta
}
}^{n}\sum_{\sigma_{i}\subseteq V_{j}}g_{i}\partial_{1}\sigma_{i}.
\end{align*}
There are $u_{\xi},u_{\eta},v\in M_{0}\left(X;G\right)$ such that
$V_{u_{\xi}}\subseteq V_{\xi}$, $V_{u_{\eta}}\subseteq V_{\eta}$,
$V_{v}\cap\left(V_{\xi}\cup V_{\eta}\right)=\varnothing$ and $\partial_{1}u=x-u_{\xi}+u_{\eta}-y+v$.
By comparing the coefficients, we have $v=0$, $u_{\xi}\neq0$ and
$u_{\eta}\neq0$. Since $\partial_{1}u=x-y$, it must hold that $u_{\xi}=u_{\eta}$.
It leads to a contradiction. Hence $x$ and $y$ are contained in
the same equivalence class of $\underset{u}{\leftrightarrow}$. There
is an infinitesimal chain connecting $x$ to $y$. By \prettyref{lem:chain-connectedness-lemma},
$X$ is chain-connected.\end{proof}
\begin{cor}
If $X$ has exactly $n$ chain-connected components, then $H_{0}\left(X;G\right)=G^{n}$.\end{cor}
\begin{proof}
The case $n=1$ is immediate from \prettyref{prop:ses-of-aug-cc}
and \prettyref{thm:0hom-of-cc-space-vanishes}. The general case can
be reduced to this case using \prettyref{prop:finite-additivity-axiom}.
\end{proof}

\section{S-homotopy and S-retraction}

We now introduce a weaker notion of homotopy. Let $X$ and $Y$ be
uniform spaces. An internal map $f:{}^{\ast}X\to{}^{\ast}Y$ is said
to be \emph{S-continuous} at $x\in{}^{\ast}X$ if $f\left(y\right)\approx f\left(x\right)$
for all $y\approx x$. For example, if a map $f:X\to Y$ is uniformly
continuous on $X$, the nonstandard extension $^{\ast}f:{}^{\ast}X\to{}^{\ast}Y$
is S-continuous on $^{\ast}X$, and vice versa. We say that two S-continuous
maps $f,g:{}^{\ast}X\to{}^{\ast}Y$ are \emph{S-homotopic} if there
is an S-continuous map $h:{}^{\ast}X\times{}^{\ast}\left[0,1\right]\to{}^{\ast}Y$,
called an \emph{S-homotopy} between $f$ and $g$, such that $h\left(\cdot,0\right)=f$
and $h\left(\cdot,1\right)=g$. If two uniformly continuous maps $f,g:X\to Y$
are uniformly homotopic, then the nonstandard extensions $^{\ast}f,{}^{\ast}g:{}^{\ast}X\to{}^{\ast}Y$
are S-homotopic. The converse is not true, e.g., the inclusion map
$i:\set{\pm1}\hookrightarrow\mathbb{R}\setminus\set{0}$ and the constant
map $c:\set{\pm1}\to\set{+1}\subseteq\mathbb{R}\setminus\set{0}$.

Every S-continuous map $f:{}^{\ast}X\to{}^{\ast}Y$ induces a homomorphism
$M_{\bullet}\left(f;G\right):M_{\bullet}\left(X;G\right)\to M_{\bullet}\left(Y;G\right)$
in the usual way. Thus the domain of $M_{\bullet}\left(\cdot;G\right)$
and $H_{\bullet}\left(\cdot;G\right)$ can be extended to the category
of uniform spaces with S-continuous maps. We obtain the following
analogue of \prettyref{prop:uniform-homotopy-axiom}.
\begin{thm}[S-homotopy]
\label{thm:S-homotopy-axiom}If two S-continuous maps $f,g:{}^{\ast}X\to{}^{\ast}Y$
are S-homotopic, then the induced homomorphisms $M_{\bullet}\left(f;G\right)$
and $M_{\bullet}\left(g;G\right)$ are chain homotopic. Hence $H_{\bullet}\left(f;G\right)=H_{\bullet}\left(g;G\right)$.
This also holds for relative $\mu$-homology.\end{thm}
\begin{proof}
Similar to \prettyref{prop:uniform-homotopy-axiom}.
\end{proof}
Let $X$ be a uniform space and $A$ a subset of $X$. Let $i:{}^{\ast}A\hookrightarrow{}^{\ast}X$
be the inclusion map. An S-continuous map $r:{}^{\ast}X\to{}^{\ast}A$
is called an \emph{S-retraction} if $ri=\mathrm{id}_{^{\ast}A}$.
An S-retraction $r:{}^{\ast}X\to{}^{\ast}A$ is called an \emph{S-deformation
retraction} if $ir$ is S-homotopic to $\mathrm{id}_{^{\ast}X}$.
\begin{thm}
\label{thm:deformation-theorem}Let $A$ be a dense subset of a uniform
space $X$. Then there is an S-deformation retraction $r:{}^{\ast}X\to{}^{\ast}A$.\end{thm}
\begin{proof}
Let $i:{}^{\ast}A\hookrightarrow{}^{\ast}X$ be the inclusion. By
saturation, there is an internal entourage $U$ of $^{\ast}X$ with
$U\subseteq\left(\approx\right)$. For any $x\in{}^{\ast}X$, by transfer,
we have $U\left[x\right]\cap{}^{\ast}A\neq\varnothing$. The family
$\set{U\left[x\right]\cap{}^{\ast}A|x\in{}^{\ast}X}$ of nonempty
sets is internal. By transfer and the axiom of choice, there is an
internal map $r:{}^{\ast}X\to{}^{\ast}A$ such that $r\left(x\right)\in U\left[x\right]\cap{}^{\ast}A$
for all $x\in{}^{\ast}X$. We can assume that $r\left(a\right)=a$
for all $a\in{}^{\ast}A$, i.e. $ri=\mathrm{id}_{^{\ast}A}$. For
any $x\approx y\in{}^{\ast}X$, since $\left(x,r\left(x\right)\right)\in U$
and $\left(y,r\left(y\right)\right)\in U$, we have that $r\left(x\right)\approx x\approx y\approx r\left(y\right)$.
Hence $r$ is S-continuous. It remains to show that $ir$ is S-homotopic
to $\mathrm{id}_{^{\ast}X}$. Fix a hyperreal number $0<\alpha<1$
and define an internal map $h:{}^{\ast}X\times{}^{\ast}\left[0,1\right]\to{}^{\ast}X$:
\[
h\left(x,t\right)=\begin{cases}
ir\left(x\right), & 0\leq t<\alpha,\\
x, & \alpha\le t\leq1.
\end{cases}
\]
Clearly $h$ is an S-homotopy between $ir$ and $\mathrm{id}_{^{\ast}X}$.
The proof is completed.
\end{proof}
The following generalizes the fact that every uniform space with a
chain-connected dense subset is chain-connected.
\begin{cor}
Let $A$ be a dense subset of a uniform space $X$. The $\mu$-homology
of $A$ is isomorphic to the $\mu$-homology of $X$.\end{cor}
\begin{proof}
Let $i:{}^{\ast}A\hookrightarrow{}^{\ast}X$ be the inclusion. By
\prettyref{thm:deformation-theorem}, there exists an S-deformation
retraction $r:{}^{\ast}X\to{}^{\ast}A$. Then $H_{\bullet}\left(r;G\right)H_{\bullet}\left(i;G\right)=H_{\bullet}\left(ri;G\right)=H_{\bullet}\left(\mathrm{id}_{^{\ast}A};G\right)=\mathrm{id}_{H_{\bullet}\left(A;G\right)}$.
By \prettyref{thm:S-homotopy-axiom}, $H_{\bullet}\left(i;G\right)H_{\bullet}\left(r;G\right)=H_{\bullet}\left(ir;G\right)=H_{\bullet}\left(\mathrm{id}_{^{\ast}X};G\right)=\mathrm{id}_{H_{\bullet}\left(X;G\right)}$.
Hence $H_{\bullet}\left(i;G\right)$ and $H_{\bullet}\left(r;G\right)$
are isomorphisms.\end{proof}
\begin{cor}
Let $\bar{X}$ be a uniform completion of a uniform space $X$. The
$\mu$-homology of $\bar{X}$ is isomorphic to the $\mu$-homology
of $X$.
\end{cor}

\begin{cor}
Let $K$ be a uniform compactification of a uniform space $X$. The
$\mu$-homology of $K$ is isomorphic to the $\mu$-homology of $X$.
\end{cor}

\section{Relation to McCord homology theory}

Recall that a member $\left(a_{0},\ldots,a_{p}\right)$ of $^{\ast}X^{p+1}$
is called a $p$-microsimplex in the sense of McCord if there is a
point $x$ of $X$ with $\set{a_{0},\ldots,a_{p}}\subseteq\mu\left(x\right)$.
This is equivalent to our definition for compact uniform spaces. Hence
$\mu$-homology theory coincides with McCord homology theory for compact
uniform spaces. This does not hold for noncompact uniform spaces.

We say that a topological space $X$ is \emph{NS-chain-connected}
if for any $x,y\in X$ there are an internal hyperfinite sequence
$\set{x_{i}}_{i=0}^{n}$ of points of $^{\ast}X$ and a sequence $\set{y_{i}}_{i=0}^{n-1}$
of points of $X$, such that $x_{0}=x$, $x_{n}=y$ and $\set{x_{i},x_{i+1}}\subseteq\mu\left(y_{i}\right)$
for all $i<n$. NS is the acronym of Near Standard. Every path-connected
space is NS-chain-connected, and every NS-chain-connected space is
connected. The converses are not true, e.g., the closed topologist's
sine curve $\set{\left(x,\sin\left(1/x\right)\right)|0<x\leq1}\cup\left(\set{0}\times\left[-1,1\right]\right)$
is NS-chain-connected but not path-connected, and the topologist's
sine curve $\set{\left(x,\sin\left(1/x\right)\right)|0<x\leq1}\cup\set{\left(0,0\right)}$
is connected but not NS-chain-connected. However, every compact connected
space is NS-chain-connected.
\begin{thm}
If a topological space is NS-chain-connected, then the $0$-th reduced
McCord homology vanishes. If the coefficient group is nontrivial,
the converse is also true.\end{thm}
\begin{proof}
Similar to \prettyref{thm:0hom-of-cc-space-vanishes} and \prettyref{thm:if-0hom-vanishes-then-cc}.\end{proof}
\begin{example}
Let $T$ be the closed topologist's sine curve. Since $T$ is NS-chain-connected,
the $0$-th reduced McCord homology of $T$ vanishes. On the other
hand, $T$ is not path-connected, so the $0$-th reduced singular
homology of $T$ does not vanish.
\end{example}
For compact uniform spaces, chain-connectedness and NS-chain-connectedness
are equivalent. There are, however, chain-connected but not NS-chain-connected
uniform spaces. Consequently, in general, $\mu$-homology theory does
not coincide with McCord's one.
\begin{example}
Let $X$ be the real line without one point. Since $X$ is chain-connected,
the $0$-th reduced $\mu$-homology of $X$ vanishes. On the other
hand, $X$ is neither connected nor NS-chain-connected, so the $0$-th
reduced McCord homology of $X$ does not vanish.
\end{example}

\section{Questions}

Let $\Delta^{p}$ be the standard $p$-simplex and $X$ a uniform
space. An S-continuous map from $^{\ast}\Delta^{p}$ to $^{\ast}X$
is called an S-singular $p$-simplex on $X$. One can define a homology
theory based on hyperfinite chains of S-singular simplices. Does this
homology theory coincide with $\mu$-homology theory?

One can also define a homotopy theory by using S-homotopy equivalence.
The higher-dimensional analogue of chain-connectedness can be formulated
in terms of this homotopy theory. Does the analogue of Hurewicz theorem
hold for this homotopy theory and $\mu$-homology theory?

\bibliographystyle{elsarticle-num}
\nocite{*}
\bibliography{\string"Nonstandard homology theory for uniform spaces\string"}

\begin{thebibliography}{1}
\expandafter\ifx\csname url\endcsname\relax
  \def\url#1{\texttt{#1}}\fi
\expandafter\ifx\csname urlprefix\endcsname\relax\def\urlprefix{URL }\fi
\expandafter\ifx\csname href\endcsname\relax
  \def\href#1#2{#2} \def\path#1{#1}\fi

\bibitem{McC72}
M.~C. McCord, Non-standard analysis and homology, Fundamenta Mathematicae
  74~(1) (1972) 21--28.

\bibitem{Gar78}
S.~Garavaglia, Homology with equationally compact coefficients, Fundamenta
  Mathematicae 100~(1) (1978) 89--95.

\bibitem{Ziv87}
R.~T. \v{Z}ivaljevi\'{c}, On a cohomology theory based on hyperfinite sums of
  microsimplexes, Pacific Journal of Mathematics 128~(1) (1987) 201--208.

\bibitem{Kor10}
T.~Korppi, On the homology of compact spaces by using non-standard methods,
  Topology and its Applications 157 (2010) 2704--2714.

\bibitem{ES45}
S.~Eilenberg, N.~E. Steenrod, Axiomatic {A}pproach to {H}omology {T}heory,
  Proceedings of the National Academy of Sciences of the United States of
  America 31~(4) (1945) 117--120.

\end{thebibliography}


\providecommand{\bysame}{\leavevmode\hbox to3em{\hrulefill}\thinspace}
\providecommand{\MR}{\relax\ifhmode\unskip\space\fi MR }
\providecommand{\MRhref}[2]{%
  \href{http://www.ams.org/mathscinet-getitem?mr=#1}{#2}
}
\providecommand{\href}[2]{#2}
\begin{thebibliography}{1}

\bibitem{Gar78}
Steven Garavaglia, \emph{Homology with equationally compact coefficients},
  Fund. Math. \textbf{100} (1978), no.~1, 89--95.

\bibitem{Ima16}
Takuma Imamura, \emph{Nonstandard homology theory for uniform spaces}, Topol.
  Appl. \textbf{209} (2016), 22--29.

\bibitem{Let95}
Steven~C. Leth, \emph{Some nonstandard methods in geometric topology},
  Developments in {N}onstandard {M}athematics (Nigel~J. Cutland, Vitor Neves,
  A.~F. Oliveira, and Jose Sousa-Pinto, eds.), Chapman and Hall/CRC, 1995,
  pp.~50--60.

\end{thebibliography}

\end{document}


\title{Corrigendum to ``Nonstandard homology theory for uniform spaces''}

\thanks{DOI of original article: \href{http://dx.doi.org/10.1016/j.topol.2016.05.016}{10.1016/j.topol.2016.05.016}.}
\thanks{The editor decided not to publish a corrigendum on the journal.}

\author{Takuma Imamura}

\address{Research Institute for Mathematical Sciences, Kyoto University, Kitashirakawa
Oiwake-cho, Sakyo-ku, Kyoto 606-8502, Japan}

\email{timamura@kurims.kyoto-u.ac.jp}

\begin{abstract}
This is an unpublished corrigendum to {\it Nonstandard homology theory for uniform spaces}, Topology and its Applications {\bf 209} (2016), 22--29. 
\end{abstract}

\maketitle
The author regrets for the errors presented in the previous paper
\cite{Ima16}, including one critical (and only partially recoverable)
error.
\begin{enumerate}
\item In Abstract and Introduction, the sentence ``every uniform space can
be S-deformation retracted to a dense subset'' is ambiguous. The correct
sentence is ``every uniform space can be S-deformation retracted to
\emph{any} dense subset''.
\item In page 23, paragraph 4 is incorrect when the index set $I$ is infinite.
In addition, there are some omissions of the star $\prescript{\ast}{}{\left({-}\right)}$.
Here is the correction: Let $X=\prod_{i\in I}X_{i}$ be a product
of uniform spaces. For every $x,y\in\prescript{\ast}{}{X}$, $x\approx_{X}y$
if and only if $x\left(i\right)\approx_{X_{i}}y\left(i\right)$ for
all $i\in I$. This equivalence is still true when the index set $I$
is infinite. Next, let $Y=\coprod_{j\in J}Y_{j}$ be a \emph{finite}
sum of uniform spaces. For every $x,y\in\prescript{\ast}{}{Y}$, $x\approx_{Y}y$
if and only if $x,y\in\prescript{\ast}{}{Y_{j}}$ and $x\approx_{Y_{j}}y$
for some $j\in J$. This equivalence is \emph{not} true when the index
set $J$ is infinite.
\item In page 24, line 5, the sentence ``$M_{\bullet}\left(f;G\right):\left(X,A;G\right)\to\left(Y,B;G\right)$''
should be ``$M_{\bullet}\left(f;G\right):M_{\bullet}\left(X,A;G\right)\to M_{\bullet}\left(Y,B;G\right)$''.
\item The proof of Proposition 5 is incomplete (and has typos). We give
a correct proof: Let $u=\sum_{j=0}^{m}g_{j}\sigma_{j}\in M_{p}\left(X;G\right)$.
Since any two points in different components of $\prescript{\ast}{}{X}$
are not infinitely close, each $\sigma_{j}$ is contained in one and
only one of the $\prescript{\ast}{}{X_{i}}$s. It follows that
\[
u=\sum_{i=0}^{n}\sum_{\substack{j=0\\
\sigma_{j}\subseteq\prescript{\ast}{}{X_{i}}
}
}^{m}g_{j}\sigma_{j}\in\sum_{i=0}^{n}M_{p}\left(X_{i};G\right)
\]
and such a decomposition is unique. Hence $M_{p}\left(X;G\right)=\bigoplus_{i=0}^{n}M_{p}\left(X_{i};G\right)$.
\item The proof of Theorem 8 is a hand-waving and has errors. For accuracy
and to correct any mistakes, the second paragraph of this proof should
be replaced with the following: Let $u=\sum_{i=0}^{n}g_{i}\sigma_{i}\in\ker\varepsilon$.
Since $\varepsilon u=0$, we have that
\[
u=\sum_{i=0}^{n}g_{i}\left(\sigma_{i}-\sigma_{0}\right)+\sum_{i=0}^{n}g_{i}\sigma_{0}=\sum_{i=0}^{n}g_{i}\left(\sigma_{i}-\sigma_{0}\right).
\]
According to the proof of Lemma 7, for each $0$-microsimplices $\sigma$
and $\tau$, we can \emph{internally} find an infinitesimal chain
connecting $\sigma$ to $\tau$. Hence for each $0\leq i\leq n$ we
can \emph{internally} find a $1$-microchain $v_{i}$ that kills $\sigma_{i}-\sigma_{0}$,
i.e., $\sigma_{i}-\sigma_{0}=\partial_{1}v_{i}$. By the internality,
the hyperfinite sum $\sum_{i=0}^{n}g_{i}v_{i}$ exists and kills $u$.
Hence $u\in\mathrm{im}\partial_{1}$.
\item In the proof of Theorem 9 (page 27, line 14), ``$\underset{u}{\leftrightarrow}$ implies $\approx$'' should be replaced with ``$x\underset{u}{\leftrightarrow}y$ implies the existence of an infinitesimal chain between $x$ and $y$''.
\item In page 29, line 13, the sentence ``the topologist's sine curve ... is not NS-chain-connected'' is certainly true but nontrivial. We denote by $\mathbb{T}$ the topologist's sine curve. Consider a moving point
$P$ on $\prescript{\ast}{}{\mathbb{T}}$ starting from $\left(1,\sin1\right)$.
$P$ can do infinitesimal jumping. To jump to $\left(0,0\right)$
from the $\sin\left(1/x\right)$-part, $P$ should infinitesimally
approach the $y$-axis, so $P$ should cross the monad of points on
the $y$-axis other than $\left(0,0\right)$. These points are remote
(i.e. not nearstandard) in $\prescript{\ast}{}{\mathbb{T}}$. Hence
$\left(1,\sin1\right)$ and $\left(0,0\right)$ cannot be connected
by a nearstandard infinitesimal chain. To justify this intuition,
we use the following lemmas:
\begin{lem}
Let $X$ be a topological space. Let $A$ be an internal subset of
$\prescript{\ast}{}{X}$. The set $\prescript{\circ}{}{A}=\set{x\in X|\mu\left(x\right)\cap A\neq\varnothing}$
is closed.
\end{lem}
\begin{proof}
Let $x$ be a closure point of $\prescript{\circ}{}{A}$. For each
open neighbourhood $U$ of $x$, there exists an $x_{U}\in\prescript{\circ}{}{A}$
such that $x_{U}\in U$. Each $x_{U}$ is infinitesimally close to
some $y_{U}\in A$. Since $U$ is open, we have that $y_{U}\in\mu\left(x_{U}\right)\subseteq\prescript{\ast}{}{U}$.
Hence $\prescript{\ast}{}{U}\cap A\neq\varnothing$ holds. By the
overspill principle, there is an internal (open) neighbourhood $U$
of $x$ such that $U\subseteq\mu\left(x\right)$ and $U\cap A\neq\varnothing$.
It follows that $\mu\left(x\right)\cap A\neq\varnothing$ and $x\in\prescript{\circ}{}{A}$.
\end{proof}
\begin{lem}
Let $X$ be a metric space. Let $A$ be a subset of $\prescript{\ast}{}{X}$.
If $A$ is S-bounded (i.e. contained in some standard ball), then
$\prescript{\circ}{}{A}$ is bounded.
\end{lem}
\begin{proof}
For some $x_{0}\in X$ and $r\in\mathbb{R}_{+}$ we have that $A\subseteq\prescript{\ast}{}{B_{r}\left(x_{0}\right)}$.
Let $x\in\prescript{\circ}{}{A}$. There exists an $a\in A$ such
that $d\left(x,a\right)\approx0$. Then, $d\left(x,x_{0}\right)\leq\prescript{\ast}{}{d}\left(x,a\right)+\prescript{\ast}{}{d}\left(a,x_{0}\right)<1+r$.
Hence $\prescript{\circ}{}{A}\subseteq B_{1+r}\left(x_{0}\right)$
and $\prescript{\circ}{}{A}$ is bounded.
\end{proof}
We here say that a subset $A$ of a topological space $X$ is \emph{connected
in the weak sense} if there are no disjoint open sets $U,V$ of $X$
such that $A\subseteq U\cup V$, $A\cap U\neq\varnothing$ and $A\cap V\neq\varnothing$.
It is obvious that every connected subset (in the usual sense) is
connected in the weak sense. The converse does not hold (e.g. the
topological space $X=\set{a,b,c}$ equipped with the topology $\set{\varnothing,\set{a,c},\set{b,c},\set{c},X}$
and the subset $A=\set{a,b}$).
\begin{lem}
Let $X$ be a topological space. Let $A$ be an internal set of nearstandard
points of $\prescript{\ast}{}{X}$. If $A$ is internally connected
in the weak sense, then $\prescript{\circ}{}{A}$ is connected in
the weak sense.
\end{lem}
\begin{proof}
Suppose not, and let $U$ and $V$ be disjoint open sets of $X$ such
that $\prescript{\circ}{}{A}\subseteq U\cup V$, $\prescript{\circ}{}{A}\cap U\neq\varnothing$
and $\prescript{\circ}{}{A}\cap V\neq\varnothing$. By transfer, $\prescript{\ast}{}{U}$
and $\prescript{\ast}{}{V}$ are disjoint internal open sets of $\prescript{\ast}{}{X}$.
For each $y\in A$, there exists an $x\in\prescript{\circ}{}{A}\left(\subseteq U\cup V\right)$
such that $y\in\mu\left(x\right)\subseteq\prescript{\ast}{}{U}\cup\prescript{\ast}{}{V}$.
Hence $A\subseteq\prescript{\ast}{}{U}\cup\prescript{\ast}{}{V}$.
Let $a\in U\cap\prescript{\circ}{}{A}$ and $b\in V\cap\prescript{\circ}{}{A}$.
There exist $a',b'\in A$ such that $a'\in\mu\left(a\right)$ and
$b'\in\mu\left(b\right)$. Since $U$ and $V$ are open, $a'\in\prescript{\ast}{}{U}$
and $b'\in\prescript{\ast}{}{V}$ hold. Hence $\prescript{\ast}{}{U}\cap A$
and $\prescript{\ast}{}{V}\cap A$ are both nonempty. This is a contradiction
and the lemma is proved.
\end{proof}
We consider $\mathbb{T}$ as a subset of the plane $\mathbb{R}^{2}$.
Now, let $\set{x_{i}}_{i=0}^{n}$ be an infinitesimal chain in $\prescript{\ast}{}{\mathbb{T}}$
that connects between $\left(1,\sin1\right)$ and $\left(0,0\right)$.
Let $A$ be the internal polygonal path $x_{0}x_{1}\cdots x_{n}$
in $\prescript{\ast}{}{\mathbb{R}}^{2}$. Then, $\prescript{\circ}{}{A}\subseteq\prescript{\circ\ast}{}{\mathbb{T}}=\mathrm{cl}\left(\mathbb{T}\right)$
holds. By the above lemmas, $\prescript{\circ}{}{A}$ is closed, bounded
and connected in the weak sense, so compact. However, no compact subset
of $\mathbb{T}$ containing $\left(1,\sin1\right)$ and $\left(0,0\right)$
is connected in the weak sense. Hence $\prescript{\circ}{}{A}$ must
have some points in $\mathrm{cl}\left(\mathbb{T}\right)\setminus\mathbb{T}$
and $\set{x_{i}}_{i=0}^{n}$ must have some remote points of $\prescript{\ast}{}{\mathbb{T}}$.
Therefore, $\mathbb{T}$ is not NS-chain-connected. See also \cite[Proposition 3 and Proposition 5]{Let95}.
\item In page 29, line 14, the sentence ``every compact connected space
is NS-chain-conneted'' is also true but nontrivial. We first prove
the uniformisable case:
\begin{prop}
Every compact, connected, uniformisable space is NS-chain-connected.
\end{prop}
\begin{proof}
Let $X$ be such a space. We equip $X$ with a (unique) uniformity
and consider $X$ as a uniform space. Then, $X$ is chain-connected.
Let $x,y\in X$. By \cite[Lemma 7]{Ima16}, there exists an infinitesimal
chain $\set{x_{i}}_{i=0}^{n}$ that connects $x$ and $y$. Since
$X$ is compact, for each $0\leq i<n$, we can find an $a_{i}\in X$
such that $x_{i}\in\mu\left(a_{i}\right)$. It follows that $x_{i+1}\in\mu\left(a_{i}\right)$.
\end{proof}
Note that this proof works in arbitrary enlargements $\prescript{\ast}{}{\mathbb{U}}$
of $\mathbb{U}$. In other words, we only need a \emph{weak} form
of the saturation principle. Next, we prove the general case:
\begin{prop}
Every compact connected space is NS-chain-connected.
\end{prop}
\begin{proof}
Since $X$ is connected, the $0$-th reduced \v{C}ech homology of
$X$ vanishes. By \cite[Theorem 2]{Gar78}, the reduced $0$-th McCord
homology of $X$ also vanishes. By the \emph{last half of} \cite[Theorem 16]{Ima16},
$X$ is NS-chain-connected.
\end{proof}
This proof works only in $\left|\mathbb{U}\right|^{+}$-saturated
elementary extensions $\prescript{\ast}{}{\mathbb{U}}$ of $\mathbb{U}$,
because Garavaglia's result depends on the \emph{full} saturation
principle (see \cite[Theorem 4]{Gar78}).
\item The \emph{first half of} Theorem 16 cannot be proved in the way of
\cite{Ima16}, and the author does not know whether it is true or false. The
same way of Theorem 8 does not work to prove the first half of Theorem
16, because we may not be able to \emph{internally} find for each
$0\leq i\leq n$ a $1$-microchain that kills $\sigma_{i}-\sigma_{0}$.
However, under the compactness, we can prove the first half of Theorem
16 as follows.
\begin{lem}
Every NS-chain-connected space is connected (page 29, line 11).
\end{lem}
\begin{proof}
Let $X$ be an NS-chain-connected space. Suppose that $X$ has a nontrivial
clopen subset $A$. Fix an $a\in A$ and a $b\in X\setminus A$. Since
$X$ is NS-chain-connected, there exists an internal hyperfinite sequence
$\set{x_{i}}_{i=0}^{n}$ in $\prescript{\ast}{}{X}$ such that $x_{0}=a$,
$x_{n}=b$ and for each $0\leq i<n$ we can find a $y_{i}\in X$ with
$\set{x_{i},x_{i+1}}\subseteq\mu\left(y_{i}\right)$. Take $i_{0}$
such that $x_{i_{0}}\in\prescript{\ast}{}{A}$ but $x_{i_{0}+1}\notin\prescript{\ast}{}{A}$.
Then, $x_{i_{0}}\in\mu\left(y_{i_{0}}\right)\cap\prescript{\ast}{}{A}\neq\varnothing$
and $x_{i_{0}+1}\in\mu\left(y_{i_{0}}\right)\cap\prescript{\ast}{}{\left(X\setminus A\right)}\neq\varnothing$.
Since $A$ and $X\setminus A$ are both closed, we see that $y_{i_{0}}\in A$
and $y_{i_{0}}\in X\setminus A$, a contradiction.
\end{proof}
\begin{prop}
If a topological space $X$ is compact and NS-chain-connected, the
reduced $0$-th McCord homology of $X$ vanishes.
\end{prop}
\begin{proof}
By the previous lemma, $X$ is connected. Hence the reduced $0$-th
\v{C}ech homology of $X$ vanishes. By \cite[Theorem 2]{Gar78},
the reduced $0$-th McCord homology of $X$ also vanishes.
\end{proof}
\end{enumerate}
The author would like to apologise for the errors and for any confusion
caused by them.

\bibliographystyle{amsplain}
\bibliography{Nonstandard homology theory for uniform spaces Corrigendum}